\documentclass[11pt]{article}
\usepackage{graphicx} 
\graphicspath{{Figures/}}
\usepackage[english]{babel}
\usepackage{spverbatim}
\usepackage{fancyvrb}
\usepackage{fvextra}
\usepackage{tabto}
\usepackage{soul}
\usepackage{xcolor}
\usepackage{amsmath}
\usepackage{amsfonts}
\usepackage{threeparttable}
\usepackage{subcaption}
\usepackage{siunitx}
\graphicspath{{figures/}}
\usepackage{textcomp}
\usepackage{lineno}
\usepackage{makecell}
\usepackage{pdfpages}
\usepackage[htt]{hyphenat}
\usepackage{amsmath}
\usepackage{amsfonts}
\usepackage{bm}
\usepackage{float}
\usepackage{bm}
\usepackage{gensymb}
\usepackage{color}
\usepackage{csquotes}
\usepackage{comment}
\usepackage{setspace}
\usepackage{titlesec}
\usepackage{authblk}
\titleformat*{\subsubsection}{\itshape}
\usepackage[normalem]{ulem}
\usepackage[colorlinks=true, linkcolor=blue, citecolor=cyan, urlcolor=black]{hyperref}
\definecolor{bg}{rgb}{0.95,0.95,0.95}

\usepackage[letterpaper,top=1in,bottom=1in,left=1.5in,right=1in,marginparwidth=1.75cm]{geometry}

\title{\Huge What is the Funniest Number? \\ An investigation of numerical humor}
\author{\Large E.~G.~Pottebaum\thanks{emily.pottebaum@yale.edu} \\ \large Department of Physics, Yale University}
\date{April 1, 2025}

\begin{document}

\maketitle

\begin{abstract}
\noindent
In a preliminary study of numerical humor, we propose the Perceived Specificity Hypothesis (PSH). The PSH states that, for nonnegative integers $<$ 100, the funniness of a number increases with its apparent precision. A survey of 68 individuals supports the veracity of this hypothesis and indicates that oddly specific numbers tend to be funniest. Our results motivate future study in this novel subfield.
\end{abstract}

\section*{Introduction}
\label{sec:intro}
This interdisciplinary study lies at the cross section of mathematics, psychology, and sociology. While some seek to find the biggest numbers \cite{bignumbers}, or define novel categories of prime numbers \cite{primes}, our aim is to describe funny numbers. What types of numbers are funny, and why? McGraw and Warren developed the Benign Violation Theory of humor, which states that humor requires some sort of violation--of a person's sense of order or supposed correctness of things--that is simultaneously perceived as benign \cite{BVT}. We use this theory to frame our examination of numerical humor.

\subsubsection*{Trivially funny numbers}
There are a few infamous numbers that are broadly considered, at least in Internet culture, to be funny. Examples of such numbers include 69 \cite{69}, 420 \cite{AP420}, and (531)8008 \cite{8008}. These are what we will call \textbf{mimetically} funny numbers, which are associated with or symbolic of a taboo or otherwise humorous topic. The humor of these numbers is not grounded in any sort of numeracy but is instead derived from the direct connection to a concrete, tangible concept. As this study is interested in the inherent humor of numbers, we take care to differentiate between numerically and mimetically funny cases. One interesting note can be made on these numbers in the Benign Violation Theory framework--numbers themselves are, despite the topic of this work, not generally thought of as funny. The fact that mimetically funny numbers are funny (insofar as they bring to mind a humorous topic) is therefore a violation of the general perception of numbers. Most people do not care enough about numbers for this violation to be anything but benign. Numbers are not supposed to be funny, so the fact that they are funny makes them funny.
\subsubsection*{Scope}
While the title of this paper suggests a search for the singular funniest number, there are simply too many numbers to determine which is the funniest of all. Instead, with this study we seek to identify patterns of humor among those numbers most recognizable to any given individual in the sampled population. From these patterns, fundamental guidelines of numerical humor may be established and potentially applied to other ranges of numbers. For this reason, we chose to limit our search to nonnegative integers with values less than 100. \\
\par
\noindent
\underline{Disclaimer}: this study uses data collected from an English-speaking and almost completely Western population, predominantly residing in the United States of America. This population of interest (POI) was selected for its accessibility to the author, but further study in other populations is strongly encouraged (see \hyperref[sec:discussion]{Discussion} for more details). Although the sampled population is usually not explicitly mentioned for sake of convenience, \textit{all generalizations and findings in this paper arise from and refer to the POI implicitly.}

\section*{Perceived Specificity Hypothesis}
\label{sec:psh}
\noindent
The colloquial, non-rigorous use of numbers in day-to-day life is often a means for providing estimations. There is an unspoken agreement to assume a level of uncertainty in many of the numbers that are used in everyday communication. We\footnote{\textit{We} here refers to people in general, not to the author specifically.} tend to round to the nearest multiple of five, at least for numbers on the order of 10, and it is reasonable to then assume a rounding uncertainty for many situations in which multiples of five are encountered. This can be demonstrated with a simple thought experiment. Imagine Person A says to Person B, ``I will be there in 10 minutes.'' If Person A arrives in 9 minutes, or 12 minutes, it is unlikely that Person B will accuse Person A of having been incorrect about their 10-minute arrival time, unless Person B is being facetiously pedantic. Now, imagine Person A instead said to Person B, ``I will be there in 11 minutes.'' 11 is not a multiple of five and carries no implicit uncertainty. Person B will assume that Person A means to arrive in exactly 11 minutes. 

\par
While 11 is numerically no more precise than 10, it certainly \textit{feels} more precise in the scenario given above. The usage of non-multiples of five violates the norm of implicit rounding in such contexts. This violation is simultaneously benign because they are just numbers and it is not that deep.\footnote{\url{https://www.urbandictionary.com/define.php?term=Not\%20that\%20deep}}
\par

\begin{figure}[H]
    \centering
    \includegraphics[scale=0.9]{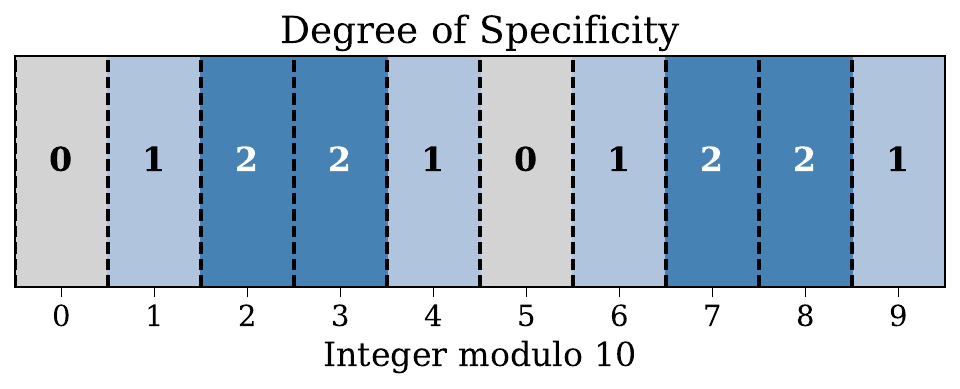}
    \caption{The degree of specificity characterizes the distance between an integer and the nearest multiple of five. Integers with a degree of specificity of 2 (dark blue) are 2 integers away from the nearest multiple of five, those with a degree of specificity of 1 (light blue) are 1 integer away from the nearest multiple of five, and those with a degree of specificity of 0 (gray) are integer multiples of five. The number line is written as the modulus of an integer with respect to 10 for generality.}
    \label{fig:deg_of_spec}
\end{figure}

We propose the Perceived Specificity Hypothesis (PSH), which states that numerical humor arises from perceived specificity. In other words, feeling more specific than is warranted is sufficient (but not necessary, see \hyperref[sec:discussion]{Discussion}) for a number to be a funny number. The concept of perceived specificity is illustrated in Figure \ref{fig:deg_of_spec}, which shows the difference between an integer (modulo 10) from the nearest multiple of five, or the \textbf{degree of specificity}.

\section*{Methodology}
\label{sec:method}
To test the Perceived Specificity Hypothesis, members of the author's professional and social circles were asked to complete a survey with several questions, all of which asked participants to select the funniest or two most funny numbers from a range of integers. The survey received a total of N = 68 responses. Most (if not all) of the respondents currently reside in North America or Europe and use English in their daily lives. Participants were informed that their answers would be used in a paper submitted to the arXiv.
\par
The survey was divided into two parts, both of which asked participants to select the funniest integer or funniest two integers from a series of different ranges, as detailed in Table \ref{tab:intranges}. In the first part, which we shall refer to as the \textbf{Zero Context} portion, participants were asked to select what they thought to be the funniest number(s). In the second part, referred to as the \textbf{Minimal Context} portion, participants were primed to consider numbers as tools of communication by being asked to consider the following scenario:
\begin{displayquote}
    Imagine you are out in public and you vaguely overhear a conversation between two strangers: \\ \\
    Stranger A asks a question, but you can't make out what the question is. Stranger B answers the question with a number. \\ \\
    Given the same options as presented previously, which response would you find the funniest? In other words, out of the provided options, which one would be the funniest answer to a question for which you have no context?
\end{displayquote}

\noindent The questions in both parts were identical except for one difference: in the Minimal Context portion, participants were required to pick only one number for all ranges instead of being able to pick two numbers from the longer ranges. The reason for this change was ostensibly to further narrow down which number was funniest, but the author acknowledges that this choice does introduce a potential design flaw that could affect comparisons between the Zero Context and Minimal Context results for the affected ranges. Additionally, there was an unintentional difference corrected partway through data collection. For 30 respondents, there are no answers for range [30, 35] in the Minimal Context portion due to a technical issue that omitted the question from the survey. This was corrected as soon as it was noticed, but the 30 people who completed the survey before the correction were not asked to retake the survey out of respect. 
\par
The selected ranges include between one and three multiples of five. To avoid decision fatigue, only one example of each range was included (e.g., instead of including [0, 5], [10, 15], [20, 25], ..., [90, 95], only the range [30, 35] was included). We strove to exclude mimetically funny numbers (e.g., 69) from the selected ranges in order to probe specifically numerical humor. 
\begin{table}[H]
    \centering
    \begin{tabular}{|c|c|}
        \hline
        \textbf{Question in Part 1 (Part 2)} & \textbf{Options} \\
        \hline
         Which number (response) is funniest? & 30, 31, 32, 33, 34, 35 \\
         \hline
         Which number (response) is funniest? & 45, 46, 47, 48, 49, 50 \\
         \hline
         Which number (response) is funniest? & 63, 64, 65, 66, 67 \\
         \hline
         \makecell[c]{Which number (response) is funniest? \\ Select up to two (select only one)} & 0, 1, 2, 3, 4, 5, 6, 7, 8, 9, 10 \\
         \hline
         \makecell[c]{Which number (response) is funniest? \\ Select up to two (select only one)} & 77, 78, 79, 80, 81, 82, 83, 84, 85, 86, 87 \\
         \hline
         
    \end{tabular}
    \caption{Abbreviated summary of the survey questions and options. ``Part 1" refers to the Zero Context portion, and ``Part 2" refers to the Minimal Context portion.}
    \label{tab:intranges}
\end{table}

\section*{Results}
\begin{figure}[H]
    \centering
    \includegraphics[width=0.93\textwidth]{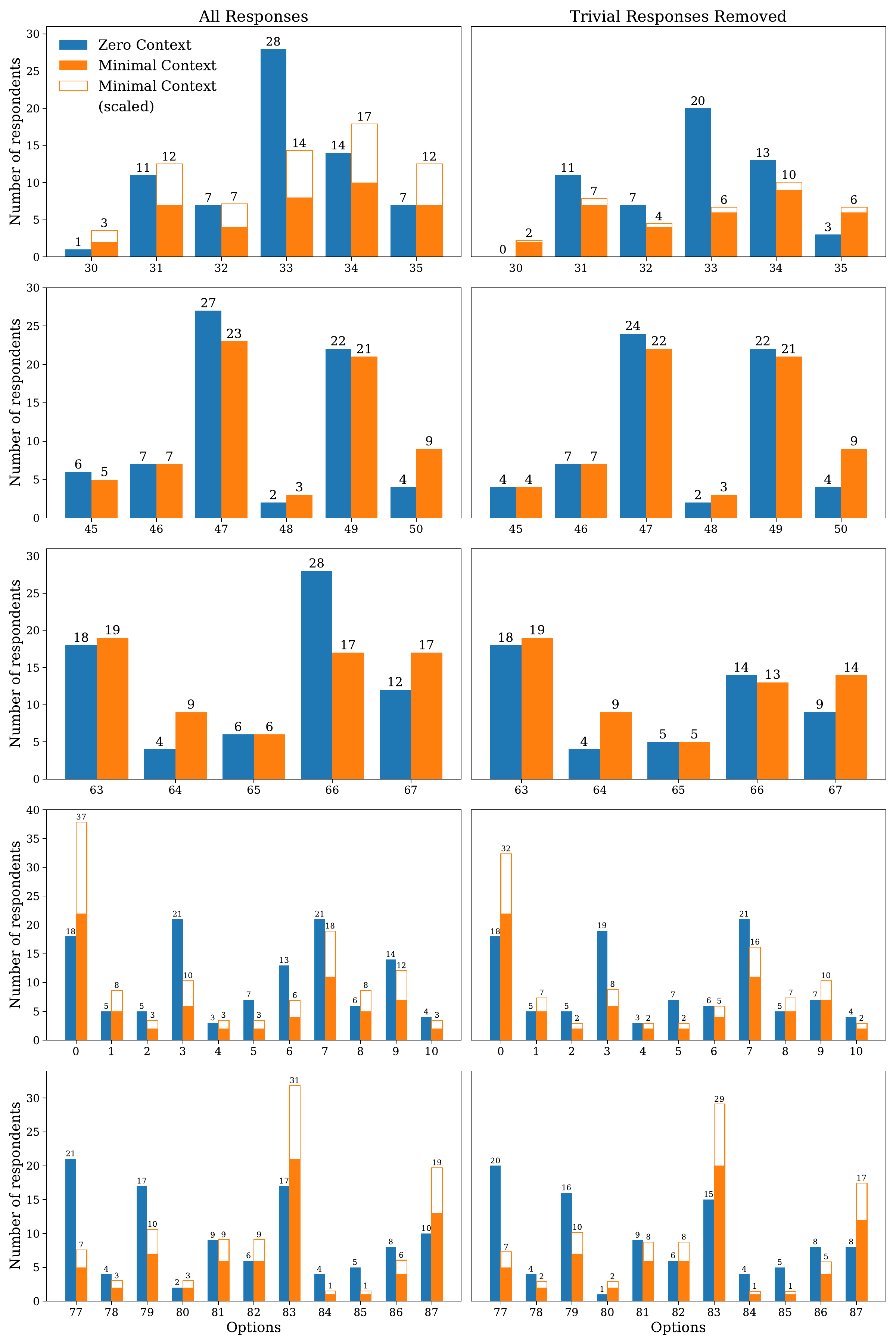}
    \caption{Survey results for ranges (from top to bottom) [30, 35], [45, 50], [63, 67], [0, 10], and [77, 87]. Results for the Zero Context portion are shown in blue, and Minimal Context results in orange (scaled to match the number of Zero Context responses when relevant). The left column shows all responses, and the right column shows the results after removing trivial responses.}
    \label{fig:results_r5}
\end{figure}

Despite efforts to exclude trivial cases, the provided explanations for several responses indicated clear mimetic motivation. The most common examples of this were explanations that relate the selected option(s) back to the numbers 69 \cite{69}, 666 \cite{666}, 25 \cite{25}, and visual association of the number 3 with breasts and/or buttocks. We thus look both at all of the responses, and at just those responses not associated with mimetically funny numbers. Because the Minimal Context portion of the survey did not include a separate optional explanation prompt for each range, just one at the end, answers were considered trivial if they matched answers deemed trivial in the Zero Context portion. When applicable, Minimal Context results were scaled proportionally with the corresponding Zero Context results (see Figure \ref{fig:results_r5}). This was done to aid in Qualitative Visual Analysis (QVA) and for relevant $\chi^2$ calculations. We do this for the ranges [30, 35], [0, 10], and [77, 87]--the first because of the 30 fewer responses than in the Zero Context portion, and the latter two ranges because participants could only select one option in the Minimal Context portion instead of two. 
\par
For full transparency, we have made the collected data publicly available with any potentially personal information redacted. The survey, data, and analysis code can be found on the author's GitHub page.\footnote{\url{https://github.com/epottebaum/funny-numbers}}

\subsection*{Observations}
\begin{figure}[h]
    \centering
    \includegraphics[width=\textwidth]{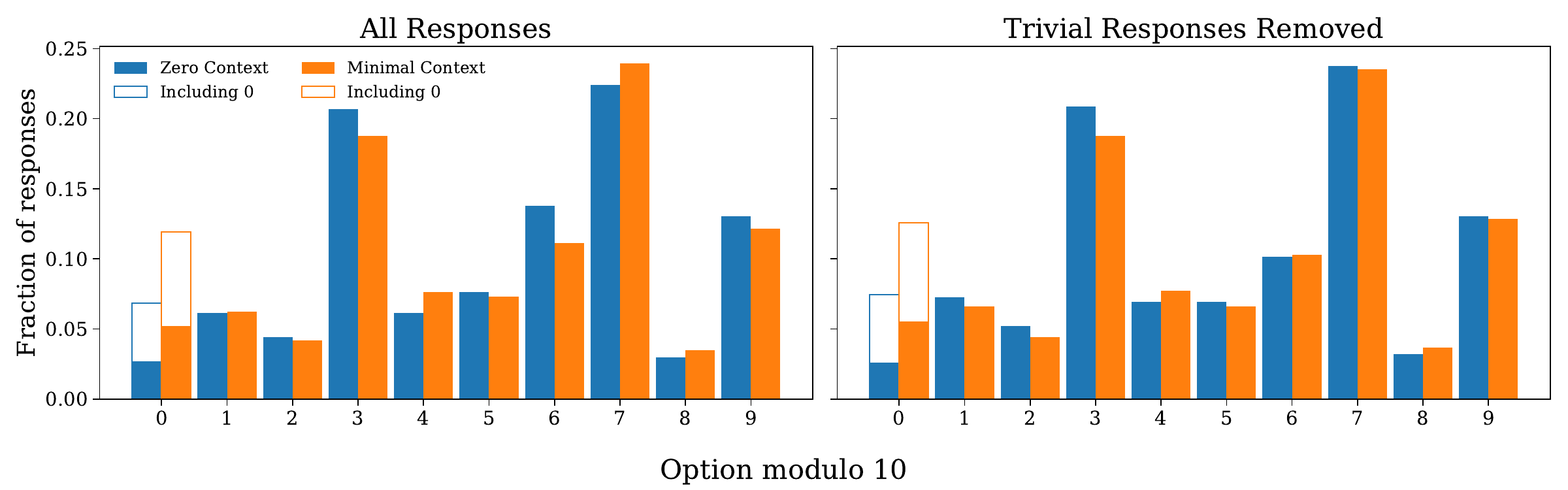}
    \caption{Combined survey results for all data (left) and data with trivial responses removed (right). Zero Context results are in blue, and Minimal Context in orange. Responses of 0 are excluded from the solid leftmost bars and included in the unfilled leftmost bars to demonstrate the popularity of 0 compared to the other options that are divisible by 10 (see \protect\hyperref[sec:zero]{\textit{Zero (0)}}).}
    \label{fig:numMod10}
\end{figure}

To study general patterns in the data, we combined survey responses for each part (i.e., for the Zero Context answers and for the Minimal Context answers) in the following way. For every integer $i$ in the given ranges (see Table \ref{tab:intranges}), we find $k$ such that $i\text{ }\%\text{ }10 = k$. The fraction of responses for each $k$ are plotted in Figure \ref{fig:numMod10}. For example, the fraction of responses for $k = 1$ includes the options 1, 31, and 81. Following are some general observations from these results.

\subsubsection*{Zero Context versus Minimal Context}
Participants showed some difference in their responses for the Zero Context portion compared to the Minimal Context portion: $\chi^2(9,N=424) = 22.20,\text{ }p = .008$ for all responses, and $\chi^2(9,N=363) = 15.57,\text{ }p = .076$ after removing trivial responses. Some participants commented on the difference between the Zero Context and Minimal Context portions in the optional prompt for explanations and general thoughts at the end of the survey:
\begin{itemize}
    \item ``Why would I change my answers"
    \item ``The pretext does not change anything"
    \item ``I’m imagining someone saying these numbers with a lot of emphasis. I can imagine someone saying `FIFTY!' pretty intensely. The eighties don’t really have the same punch to them."
    \item ``I have no idea what you're talking about"
    \item ``My first thoughts upon beginning to answer these questions were: `What is the question?' I devised many possible answers which were actually questions. Which led me to think to myself: like Jeopardy! Of course following that train of thought led me to find the funniest answers out of context, while these may be different from the previous responses I've given, I am absolutely positive that they are each individually the funniest and most beautiful numbers in existence. You do not get to decide which number is the funniest in each moment, however you do have the power to channel the comedy of mathematics through you in every second. The cosmic random number generator that courses through each of our beings determines what number will be funniest, and changes its mind frequently. I think the funniest number is neither eternal nor fleeting, it simply is, was, will be again, and never has been. That is the best answer I can give, and I thank you for your careful consideration"
    \item ``I don't really know sorry."
\end{itemize}
\noindent
Because the number of responses for some questions was not consistent between the Zero and Minimal Context parts of the survey (see \hyperref[sec:method]{Methodology}), the Minimal Context data were scaled in proportion to the number of Zero Context responses for the $\chi^2$ calculations. We thus refrain from drawing conclusions from the direct comparison between the Zero and Minimal Context to avoid sacrificing statistical rigour.

\subsubsection*{Zero (0)}
\label{sec:zero}
Our data indicate that, when compared to the other options that are divisible by 10 (10, 30, 50, 80), participants found 0 to be significantly funnier: $\chi_{zc}^2(1, N=29) = 10.71,\text{ }p = .001$ and $\chi_{mc}^2(1, N=37) = 12.93,\text{ }p < .001$, where $zc$ refers to Zero Context data and $mc$ refers to Minimal Context data. This significance is strictly larger for data with trivial responses removed, and in both cases for Minimal Context compared to Zero Context. For comparison, any particular preferences in the funniness distribution (i.e., number of responses) between 10, 30, 50, and 80 were much less significant: $\chi_{zc}^2(3, N=11) = 2.45,\text{ }p = .484$ and $\chi_{mc}^2(3, N=15) = 9.80,\text{ }p = .020$. QVA of Figure \ref{fig:results_r5} indicates a slight preference for 50 over 30 and 80 in the Minimal Context data.
\par
A few different explanations (reworded for clarity) were provided by some who selected 0 in the Minimal Context portion: 
\begin{itemize}
    \item It is funnier to use the word ``zero" than the word ``none" in a conversation
    \item In response to a question, the smallest or largest number implies an extreme [and, presumably, a more extreme response is funnier than a less extreme response]
    \item A number is funnier in response to a question [in the Minimal Context scenario] if it isn't a common answer to an everyday question. For example, if the question was about ranking something, 0 is funnier than the rest of the integers in [1, 10] because it indicates someone felt the need to leave the usual scale of 1 to 10 to express their pronounced dislike or disapproval.
\end{itemize}

\subsubsection*{Even and odd numbers}
According to the survey, odd numbers are consistently funnier than even numbers: $\chi_{zc}^2(1, N=424) = 79.25,\text{ }p < .001$ and $\chi_{mc}^2(1, N=310) = 40.93,\text{ }p < .001$ for all data, and $\chi_{zc*}^2(1, N=363) = 76.54,\text{ }p < .001$ and $\chi_{mc*}^2(1, N=294) = 37.62,\text{ }p < .001$ for data with trivial responses removed (indicated by $*$). Some of the provided explanations that mention even and/or odd numbers include:
\begin{itemize}
    \item ``Even numbers are boring''
    \item ``For a number to be funny, it must be odd''
        \begin{itemize}
            \item The same participant later gave a related explanation: ``0 is objectively the funniest number because it is even, therefore breaking a previous rule." This is a clear and fascinating example of the Benign Violation Theory's applicability to numerical humor.
        \end{itemize}
    \item ``Evens are very to-the-point" [therefore serious, therefore not funny]
    \item Numbers with both even and odd digits are funny [e.g., 83]
    \item ``Odd numbers go brrrrr"\footnote{Likely a reference to the \textit{haha X go brr} Internet meme, see \url{https://knowyourmeme.com/memes/money-printer-go-brrr}}
    \item ``It's odd" [therefore funny]
\end{itemize}

\subsection*{Testing the Perceived Specificity Hypothesis}
\begin{figure}[H]
    \centering
    \includegraphics[width=\textwidth]{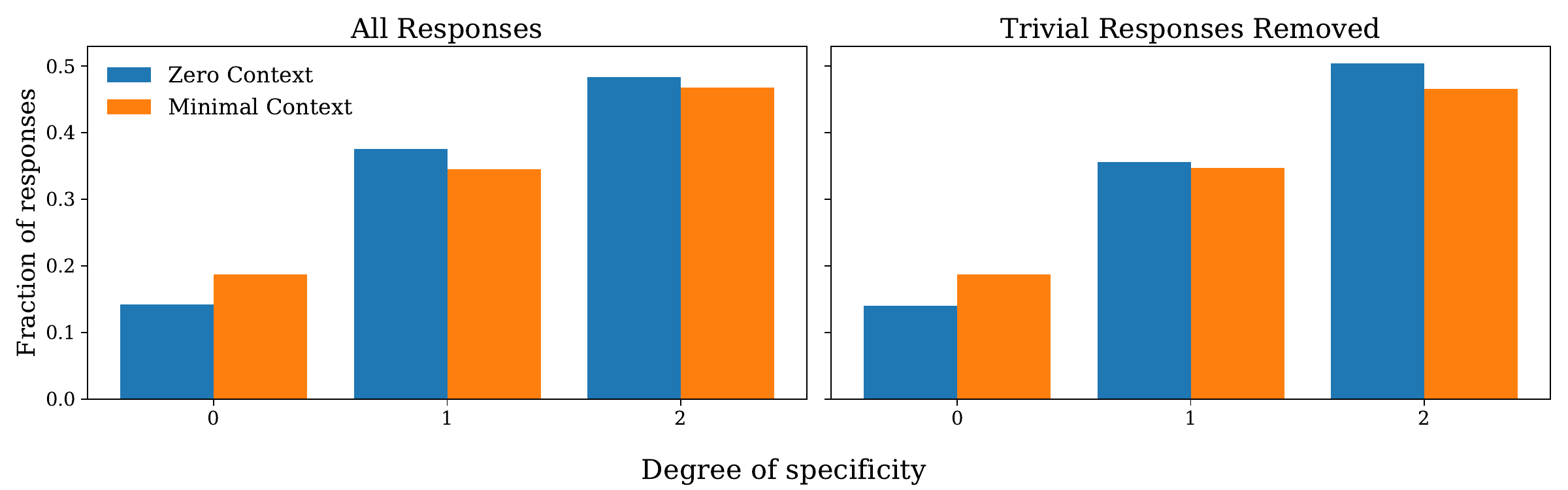}
    \caption{The fraction of responses according to the selected option's degree of specificity for all responses (left) and not including trivial responses (right). Zero Context data is shown in blue, Minimal Context data in orange.}
    \label{fig:results_genDos}
\end{figure}

To test the PSH, we consider the fraction of responses for the combined survey data (plotted in Figure \ref{fig:numMod10}) according to the degree of specificity \hyperref[fig:deg_of_spec]{as defined previously}. According to the survey, the funniness of a number is associated with the degree of specificity: $\chi_{zc}^2(2, N=424) = 15.30,\text{ }p < .001$ and $\chi_{mc}^2(2, N=310) = 6.15,\text{ }p = .046$. This association is less significant in the Minimal Context data, although the significance substantially increases for both sets of data upon removing responses of 0: $\chi_{mc,no\text{ }0}^2(2,N = 288) = 16.39,\text{ }p < .001$. The level of significance is consistent after removing trivial responses. It is clear from Figure \ref{fig:results_genDos} that participants were more likely to select numbers in the given range with higher degrees of specificity.
\par

\noindent
We also find a difference between even and odd numbers in the relationship between a number's funniness and its degree of specificity. Assigning odd numbers to negative degrees of specificity and even numbers to positive degrees of specificity, participants were far more likely to find a number with a degree of specificity of $\pm$2 funny if the number is odd: 
$$\chi_{zc}^2(1, N=205) = 102.56,\text{ }p = 4.18 \times 10^{-24}\text{ }(< .001)$$ 
$$\chi_{mc}^2(1, N=145) = 70.35,\text{ }p = 4.96 \times 10^{-17}\text{ }(< .001)$$
and, after removing trivial responses, 
$$\chi_{zc*}^2(1, N=183)=85.38,\text{ }p = 2.46 \times 10^{-20}\text{ }(< .001)$$ 
$$\chi_{mc*}^2(1, N=137)=63.13,\text{ }p = 1.93 \times 10^{-15}\text{ }(< .001)$$
\noindent
Interestingly, survey results show no preference between even and odd numbers with a degree of specificity of $\pm$1:
$$\chi_{zc}^2(1, N=159) = 0.057,\text{ }p = .812$$ 
$$\chi_{mc}^2(1, N=107) = 0.0093,\text{ }p = .923$$
$$\chi_{zc*}^2(1, N=129)=0.94,\text{ }p = .333$$ 
$$\chi_{mc*}^2(1, N=102)=0.16,\text{ }p = .692$$
\begin{figure}[H]
    \centering
    \includegraphics[width=\textwidth]{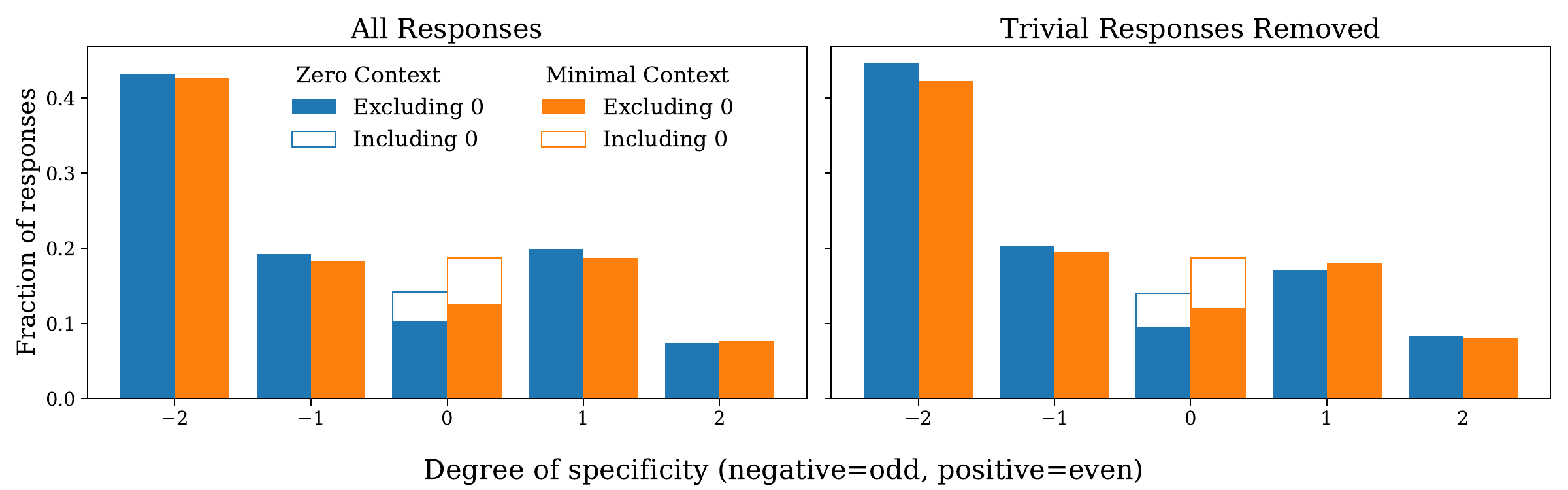}
    \caption{The fraction of responses according to the selected option's degree of specificity, which here is negative for odd numbers and positive for even numbers. Results are shown with responses of 0 both included (unfilled bars) and excluded (solid bars).}
    \label{fig:degofspecResults}
\end{figure}
\noindent
Some relevant explanations and comments from participants:
\begin{itemize}
    \item ``[Funnier numbers are those that] seem like very specific numbers rather than rounded to an even number or multiple of five"
    \item ``Numbers that are not multiples of 10, 5, and 2 are funniest."
    \item ``I think the thing that makes the number funny out of context is the entropy of the number, how specific/precise it is. 100 is a very unfunny and serious reply."
    \item ``49 is funny because it's not 50, [...] it's specifically 49!"
    \item ``The numbers that feel randomer feel funnier"
    \item ``Oddly specific numbers are always funnier"
\end{itemize}

\section*{Discussion}
\label{sec:discussion}
The results of this study support the PSH, suggesting that people in the POI tend to find that numbers further from multiples of five are funniest, at least for the most commonly used integers. However, we would like to reiterate that this is by no means an exhaustive study. While the PSH seems to be successful in explaining why some numbers are funny, we have shown evidence that at least one number, 0, is both non-mimetically funny and a multiple of five--it follows then that there must be at least one other explanation for non-mimetic numerical humor. How might the humor of 0 be explained? Could this explanation apply to other numbers? Additionally, what patterns of humor can be found in non-integer numbers? The PSH describes numbers that ``feel" more precise; is there a corresponding correlation between funniness and mathematical precision?
\par
These are just a few of many potential considerations for future study. Several survey respondents explained their choices of the funniest numbers based on the visual appearance of the number. It would be interesting to see results of the same survey written in braille, translated into languages other than English, and given verbally or in sign language. Most of these options would also necessarily study numerical humor in populations other than the POI. Research shows that the way humor tends to be used and perceived can vary between cohorts \cite{cultural-humor,cultural-humor2,cultural-humor3,cultural-humor4}, so it cannot be assumed that standards of numerical humor are universal.

\section*{Conclusion}
We found that, for nonnegative integers $<$ 100, individuals in the population of interest tend to find the numbers with a higher degree of specificity to be funnier. This can be explained by the Perceived Specificity Hypothesis, which states that apparent precision begets numerical humor. Furthermore, the most oddly specific numbers--odd numbers with a degree of specificity of 2--are the most funny, according to the data presented here. The Perceived Specificity Hypothesis is itself a specific case of the Benign Violation Theory \cite{BVT} applied to everyday use of numbers in the population of interest. While this study does not claim to have discovered the funniest number, nor a universal root of numerical humor, we have developed an initial theory of funny numbers, and hope to see future study in this field.

\section*{Acknowledgements}
I would like to thank Isaac and Ryan for introducing me to the topic of big numbers, and to Ryan and Claire for engaging in initial discussion on funny numbers. Thanks also to Andrew and Tom for inadvertently beta testing my data collection procedure with their valuable insights in our following conversations. I acknowledge my Ph.D. advisor, who I shall not name out of respect for her academic integrity, for her exasperation upon learning about this study. I thank her for putting up with my antics and plead that she continue to do so until I graduate. Finally, a monumental thank you to all who participated in the funny numbers survey, without whom this study would not have been possible. 

\raggedright
\bibliographystyle{unsrt}
\bibliography{references}

\begin{thebibliography}{10}

\bibitem{bignumbers}
Big number duel.
\newblock Googology Wiki.
\newblock \url{https://googology.fandom.com/wiki/Big_Number_Duel}.

\bibitem{primes}
Category:classes of prime numbers.
\newblock Wikipedia.
\newblock \url{https://en.wikipedia.org/wiki/Category:Classes_of_prime_numbers}.

\bibitem{BVT}
A.P. McGraw and C.~Warren.
\newblock {Benign violations: Making immoral behavior funny}.
\newblock {\em Psychological Science}, 21(8):1141--1149, 2010.

\bibitem{69}
Brian Feldman.
\newblock {Why 69 is the Internet's Coolest Number (Sex)}.
\newblock {\em Intelligencer}, June 2016.
\newblock \url{https://nymag.com/intelligencer/2016/06/why-69-is-the-internets-coolest-number-sex.html}.

\bibitem{AP420}
Gene Johnson.
\newblock How 4/20 grew from humble roots to marijuana's high holiday.
\newblock {\em Associated Press}, April 2024.
\newblock \url{https://apnews.com/article/420-history-marijuana-07ecdd87c4756cc4a1dbd4ad3b1f682c}.

\bibitem{8008}
Tom Dalzell and Terry Victor.
\newblock {\em {The Concise New Partridge Dictionary of Slang and Unconventional English}}.
\newblock Taylor \& Francis, November 2014. p.\,2060.
\newblock \url{https://books.google.com/books?id=Ak6cBQAAQBAJ&q=5318008&pg=PA2060#v=snippet&q=5318008&f=false}.

\bibitem{666}
Roland Martin.
\newblock Number of the beast.
\newblock Britannica.
\newblock \url{https://www.britannica.com/topic/number-of-the-beast}.

\bibitem{25}
YK~Eagle.
\newblock {SpongeBob funnier than 24 (Hillarious) [sic]}.
\newblock YouTube.
\newblock [Video] \url{https://www.youtube.com/watch?v=xfFlVh1IZFg&ab_channel=YKEagle}.

\bibitem{cultural-humor}
Jiang T., Li~H., and Hou Y.
\newblock {Cultural Differences in Humor Perception, Usage, and Implications}.
\newblock {\em Front Psychol}, 10(123), Jan 2019.

\bibitem{cultural-humor2}
Yue Xiaodong, Jiang Feng, Lu~Su, and Hiranandani Neelam.
\newblock {To Be or Not To Be Humorous?} {Cross Cultural Perspectives on Humor}.
\newblock {\em Frontiers in Psychology}, 7, Oct 2016.

\bibitem{cultural-humor3}
Gil Greengross.
\newblock {Humor and Aging - A Mini-Review}.
\newblock {\em Gerontology}, 59(5):448--453, Aug 2013.

\bibitem{cultural-humor4}
Darren~David Chadwick and Tracey Platt.
\newblock {Investigating Humor in Social Interaction in People With Intellectual Disabilities: A Systematic Review of the Literature}.
\newblock {\em Frontiers in Psychology}, 9, Sep 2018.

\end{thebibliography}
\end{document}